\documentclass[12pt,a4paper]{tac}
\usepackage{amssymb, amsmath, enumerate, microtype, stmaryrd, mathrsfs}
\usepackage[latin1]{inputenc}
\usepackage[arrow, matrix, tips, curve]{xy}
\SelectTips{cm}{10}

 \DeclareMathOperator{\colim}{colim}
 
\DeclareMathOperator{\el}{el} 

\newcommand{\cat}[1]{\mathbf{#1}} 

\newcommand{\op}{\mathrm{op}} 

\newcommand{\thg}{{\mathord{\text{--}}}}

\newcommand{\abs}[1]{{\left|{#1}\right|}}

\newcommand{\cd}[2][]{\vcenter{\hbox{\xymatrix#1{#2}}}}
\newcommand{\cdl}[2][]{\xymatrix@1#1{#2}}

\newcommand{\A}{{\mathcal A}}
\newcommand{\B}{{\mathcal B}}
\newcommand{\C}{{\mathcal C}}
\newcommand{\D}{{\mathcal D}}
\newcommand{\E}{{\mathcal E}}
\newcommand{\F}{{\mathcal F}}

\newcommand{\I}{{\mathcal I}}
\newcommand{\J}{{\mathcal J}}

\newcommand{\M}{{\mathcal M}}

\renewcommand{\P}{{\mathcal P}}

\newcommand{\Ss}{{\mathcal S}}
\newcommand{\T}{{\mathcal T}}

\newcommand{\Env}[1][\C]{\mathcal I#1}
\newcommand{\En}{\mathcal I}

\newcommand{\xtor}[1]{\cdl[@1]{{} \ar[r]|-{\object@{|}}^{#1} & {}}}

\newcommand{\twocong}[2][0.5]{\ar@{}[#2] \save ?(#1)*{\cong}\restore}
\newcommand{\rtwocell}[3][0.5]{\ar@{}[#2] \ar@{=>}?(#1)+/l 0.2cm/;?(#1)+/r 0.2cm/^{#3}}
\newcommand{\rtwocello}[3][0.5]{\ar@{}[#2] \ar@{=>}?(#1)+/l 0.2cm/;?(#1)+/r 0.2cm/_{#3}}
\newcommand{\ltwocell}[3][0.5]{\ar@{}[#2] \ar@{=>}?(#1)+/r 0.2cm/;?(#1)+/l 0.2cm/^{#3}}
\newcommand{\ltwocello}[3][0.5]{\ar@{}[#2] \ar@{=>}?(#1)+/r 0.2cm/;?(#1)+/l 0.2cm/_{#3}}
\newcommand{\dtwocell}[3][0.5]{\ar@{}[#2] \ar@{=>}?(#1)+/u  0.2cm/;?(#1)+/d 0.2cm/^{#3}}
\newcommand{\dthreecell}[3][0.5]{\ar@{}[#2] \ar@3{->}?(#1)+/u  0.2cm/;?(#1)+/d 0.2cm/^{#3}}
\newcommand{\utwocell}[3][0.5]{\ar@{}[#2] \ar@{=>}?(#1)+/d 0.2cm/;?(#1)+/u 0.2cm/_{#3}}
\newcommand{\dtwocelltarg}[3][0.5]{\ar@{}#2 \ar@{=>}?(#1)+/u  0.2cm/;?(#1)+/d 0.2cm/^{#3}}
\newcommand{\utwocelltarg}[3][0.5]{\ar@{}#2 \ar@{=>}?(#1)+/d  0.2cm/;?(#1)+/u 0.2cm/_{#3}}

\newdir{ >}{{}*!/-5pt/\dir{>}}

\makeatletter

\renewcommand{\paragraph}{\@startsection
{paragraph}%
{3}%
{0mm}%
{-\baselineskip}%
{-0.4em plus 0.2em minus 0.2em}%
{\normalfont\normalsize\bfseries}}%
\makeatother

\numberwithin{equation}{section} \numberwithin{paragraph}{section}

\newenvironment{Thm}{\paragraph{Theorem:}\em}{\vskip\baselineskip}
\newenvironment{Defn}{\paragraph{Definition:}}{\vskip\baselineskip}

\newenvironment{Prop}{\paragraph{Proposition:}\em}{\vskip\baselineskip}

\newenvironment{Lemma}{\paragraph{Lemma:}\em}{\vskip\baselineskip}
\newenvironment{Exs}{\paragraph{Examples:}}{\vskip\baselineskip}

\makeatletter \@namedef{itemize*}{\itemize\parsep\z@ \parskip\z@}
\@namedef{enditemize*}{\enditemize} \@namedef{enumerate*}{\enumerate\parsep\z@
\parskip\z@} \@namedef{endenumerate*}{\endenumerate}
\def\Pr@@f{\subsubsection*{\textbf{Proof}}}
\def\pr@@f[#1]{\subsubsection*{{\textbf{Proof}} #1}}
\makeatother
\DeclareMathAlphabet      {\mathbf}{OT1}{cmr}{b}{n}

\begin{document}

\copyrightyear{2014}

\title{The Isbell monad} \author{Richard Garner}
\eaddress{richard.garner@mq.edu.au} \date{\today}  \thanks{The support of an Australian
  Research Council Discovery Project, grant number DP110102360, is gratefully acknowledged.}
\maketitle
\begin{abstract}
  In 1966~\cite{Isbell1966Structure}, John Isbell introduced a
  construction on categories which he termed the ``couple category''
  but which has since come to be known as the \emph{Isbell envelope}.
  The Isbell envelope, which combines the ideas of contravariant and
  covariant presheaves, has found applications in category theory,
  logic, and differential geometry. We clarify its meaning by
  exhibiting the assignation sending a locally small category to its
  Isbell envelope as the action on objects of a pseudomonad on the
  $2$-category of locally small categories; this is the \emph{Isbell
    monad} of the title. We characterise the pseudoalgebras of the
  Isbell monad as categories equipped with a \emph{cylinder
    factorisation system}; this notion, which appears to be new, is an
  extension of Freyd and Kelly's notion of factorisation
  system~\cite{Freyd1972Categories} from orthogonal classes of arrows
  to orthogonal classes of cocones and cones.
\end{abstract}

\section{Introduction}
One of the most fundamental constructions in category theory is that
which assigns to a small category $\C$ the Yoneda embedding $Y \colon
\C \to [\C^\op, \cat{Set}]$ into its category of presheaves. As is
well known, this embedding has the effect of exhibiting $[\C^\op,
\cat{Set}]$ as a \emph{free cocompletion} of $\C$: the value at $\C$
of a left biadjoint
\begin{equation}
\cd[@C+1.5em]{
\cat{COCTS} \ar@<-1pt>[r]_{} \ar@<4.25pt>@{}[r]|{\bot} \ar@<0pt>@/^0.8em/@{<-}[r]^{} &
\cat{CAT  } \ \ \ 
}\label{eq:1}
\end{equation}
to the forgetful $2$-functor from small-cocomplete categories and
cocontinuous functors to locally small ones. At a $\C$ which is not
necessarily small, this left biadjoint still exists, but now has its
unit $Y \colon \C \to \P \C$ given by the Yoneda embedding into the
subcategory $\P \C \subset [\C^\op, \cat{Set}]$ of \emph{small}
presheaves: those which can be expressed as small colimits of
representables. Composing the two biadjoints in~\eqref{eq:1} exhibits
the process of free cocompletion as the functor part of a pseudomonad
$\P$ on $\cat{CAT}$, and it turns out that the $\P$-pseudoalgebras and
algebra pseudomorphisms are once again the small-cocomplete categories
and cocontinuous functors between them; 
which is to say that the biadjunction~\eqref{eq:1} is
\emph{pseudomonadic}~\cite{Le-Creurer2002Becks}.



Dually, we speak of \emph{free completions} of categories, meaning the
values of a left biadjoint to the forgetful $2$-functor $\cat{CTS} \to
\cat{CAT}$ from complete categories to locally small ones. The free
completion of a small $\C$ is witnessed by the dual Yoneda embedding
$Y \colon \C \to [\C, \cat{Set}]^\op$, while the general completion $Y
\colon \C \to \P^\dagger \C$ is  constructed as $\P^\dagger \C =
\P(\C^\op)^\op \subset [\C, \cat{Set}]^\op$. As before, the
biadjunction $\cat{CTS} \leftrightarrows \cat{CAT}$ induced by free
completion is pseudomonadic, so that, as before, complete categories
and continuous functors between them may be identified with
$\P^\dagger$-pseudoalgebras and their pseudomorphisms.

In~\cite[\S1.1]{Isbell1966Structure}, Isbell describes a construction
that, in some sense, combines the processes of free completion and
cocompletion; while Isbell calls this construction the ``couple
category'', we follow Lawvere in terming it the \emph{Isbell
  envelope}. Given a locally small category $\C$, the objects of its Isbell envelope
$\Env$ are
  triples $(X^+, X^-, \xi^X)$ where $X^+ \in \P\C$ and
  $X^- \in \P^\dagger\C$ and
$\xi^X_{ab} \colon X^-(b) \times X^+(a) \to \C(a,b)$
is a family of functions, natural in $a$ and $b$; while morphisms
  $(X^+, X^-, \xi^X) \to (Y^+, Y^-, \xi^Y)$ in $\Env$ are
  pairs $(f^+, f^-)$, where $f^+ \colon X^+ \to Y^+$ in $\P\C$ and
  $f^- \colon X^- \to Y^-$ in $\P^\dagger \C$ are such that each square
  \begin{equation}
\cd[@C+0.2em]{
Y^-(b) \times X^+(a) \ar[r]^-{1 \times f^+} \ar[d]_{f^- \times 1} & 
Y^-(b) \times Y^+(a) \ar[d]^{\xi^Y}\\
X^-(b) \times X^+(a) \ar[r]^-{\xi^X} & \C(a,b)
}\label{eq:3}
\end{equation}
commutes in $\cat{Set}$. There is a Yoneda embedding $Y \colon \C \to
\Env$ into the Isbell envelope, whose value at an object $c$ is given
by:
\[
\big( \ \C(\thg, c) \in [\C^\op, \cat{Set}],\ \ \C(c, \thg) \in [\C,
\cat{Set}],\ \  (\C(c,b) \times \C(a,c) \xrightarrow{\circ}
\C(a,b))_{a,b} \ \big)\rlap{ ,}
\]
and it is related to the usual two Yoneda embeddings of $\C$ through
projection functors $\pi_1$ and $\pi_2$ fitting into a commuting diagram
\begin{equation}
\cd[@+0.3em@C+1em]{
& \C \ar[d]^Y \ar[dl]_Y \ar[dr]^Y \\ \P\C & \Env \ar[l]_-{\pi_1}
\ar[r]^-{\pi_2} & \P^\dagger\C\rlap{ .}}\label{eq:4}
\end{equation}

Isbell envelopes have a range of applications. Isbell used them to
study \emph{normal completions} of categories~\cite{Isbell1967Normal}
(the categorical correlate of Dedekind--MacNeille completions of
posets); they are closely related to constructions in linear
logic~\cite{Barr1979ast-autonomous,Pratt2010Communes}, due in part to
the ``self-duality'' $\Env \cong \Env[(\C^\op)]^\op$;
in~\cite{Stacey2011Comparative} they were used to study convenient
categories of smooth spaces; and in future work we will see that they
play a role in the \emph{Reedy categories}~\cite{Reedy1974Homotopy} of
abstract homotopy theory\footnote{Roughly speaking, if $\C$ is a Reedy
  category, then an element of the Isbell envelope $\Env$ is what one
  needs to extend $\C$ to a Reedy category with one additional
  object.}. 
In this paper, however, our interest in Isbell envelopes stems from
the following natural question: given that the two outside Yoneda
embeddings in~\eqref{eq:4} are the units at $\C$ of the pseudomonads
for small-cocomplete and small-complete categories, is there a
corresponding pseudomonad whose unit at $\C$ is the central embedding?
The main contribution of this paper is answer this question in the
affirmative; the pseudomonad in question is the \emph{Isbell monad} of
the title, and we will characterise its pseudoalgebras as categories
equipped with a \emph{cylinder factorisation system}.

By a \emph{cylinder} between small diagrams $D \colon \I \to \C$ and
$E \colon \J \to \C$, we mean a family of maps $r = (r_{ij} \colon Di
\to Ej)$ natural in $i$ and $j$. A cylinder factorisation system
provides a way of factorising each such cylinder in an
essentially-unique way as a cocone followed by a cone; the unicity is
assured by the requirement that the two parts of the factorisation
should lie in suitably orthogonal classes $\E$ of cocones and $\M$ of
cones. Cylinder factorisation systems are thus a generalisation of the
\emph{orthogonal factorisation systems} of~\cite{Freyd1972Categories}
from single maps to small families of maps; while certain aspects of
this generalisation are known in the literature, the complete
definition appears to be new; we give it in
Section~\ref{sec:cylind-fact-syst}.

Now our first main result, Theorem~\ref{thm:mainresult},
exhibits a biadjunction
\begin{equation}
\cd[@C+1.5em]{
\cat{CFS} \ar@<-1pt>[r]_{} \ar@<4.5pt>@{}[r]|{\bot} \ar@<0pt>@/^0.8em/@{<-}[r]^{} &
\cat{CAT}
}\label{eq:15}
\end{equation}
between categories and cylinder factorisation systems on categories,
with as unit at $\C$ the embedding $Y \colon \C \to \Env$
of~\eqref{eq:4}. Composing the biadjoints, we thus exhibit this
embedding as the unit at $\C$ of a pseudomonad on $\cat{CAT}$, which
is the Isbell monad we seek. Our second main result,
Theorem~\ref{thm:pseudo}, shows that the pseudoalgebras for the Isbell
monad correspond with categories equipped with cylinder factorisation
systems; in other words, we show that~\eqref{eq:15},
like~\eqref{eq:1}, is pseudomonadic. This
generalises~\cite{Korostenski1993Factorization}'s characterisation of
orthogonal factorisation systems as pseudoalgebras for the
squaring monad $(\thg)^{\cat 2}$ on $\cat{CAT}$.

Our third main result concerns morphisms of cylinder factorisation systems,
of which we have said nothing so far. Given categories $\C$ and $\D$
equipped with cylinder factorisation systems, the morphisms between
them in $\cat{CFS}$ are functors $F \colon \C \to \D$ preserving both
the $\E$-cocones and the $\M$-cones; part of the pseudomonadicity
result is that these correspond with the pseudomorphisms of Isbell
pseudoalgebras. However, we also have the more general notion of
\emph{lax} and \emph{colax} morphisms of pseudoalgebras; and
Theorem~\ref{thm:3} shows that these correspond to functors $F \colon
\C \to \D$ preserving \emph{only} $\M$-cones or $\E$-cocones
respectively.

We conclude the paper by discussing variants of the notion of cylinder
factorisation systems involving factorisations for only \emph{certain}
kinds of cylinders; our final main result, Theorem~\ref{thm:4},
exhibits these as the pseudoalgebras for certain variants of the
Isbell monad, obtained by constraining the presheaves $X^+ \in \P\C$
and $X^- \in \P^\dagger \C$ that constitute an object of $\Env$ to lie
in suitable \emph{saturated classes}~\cite{Albert1988The-closure} of
weights for colimits and limits.

\section{Cylinder factorisation systems}
\label{sec:cylind-fact-syst}

Suppose that $D \colon \I \to \C$ and $E \colon \J \to \C$ are
diagrams in a category $\C$. By a \emph{cocone under $D$ with vertex
  $V$}, we mean a natural transformation $p \colon D \to \Delta V$
into the constant functor at $V$, and by a \emph{cone over $E$ with
  vertex $W$}, a natural transformation $q \colon \Delta W \to E$.
Given a map $f \colon V \to W$, we may postcompose $p$ or precompose
$q$ with it to obtain a cocone $f \cdot p \colon D \rightarrow \Delta
W$ or cone $q \cdot f \colon \Delta V \to E$.
By a \emph{cylinder from $D$ to $E$}, written $r \colon D
\rightsquigarrow E$, we mean a natural transformation
\[
\cd[@!@-2.4em@R-1.3em]{
& \I \ar[dr]^D 
\dtwocell{dd}{r} \\ \I \times \J \ar[ur]^-{\pi_1} \ar[dr]_-{\pi_2} & &
\C\rlap{ ,}
\\
& \J \ar[ur]_E
}
\]
thus, a natural family of maps $(r_{ij} \colon Di \to Ej)_{i,j \in \I
  \times \J}$. For example, if $\J = 1$, then $E$ picks out a single
vertex and so a cylinder is simply a cocone; while if $\I=1$ then a
cylinder is just a cone. For a further example, if $p \colon D
\to \Delta V$ is a cocone and $q \colon \Delta V \to
E$ a cone, then we have a cylinder $q \cdot p \colon D
\rightsquigarrow E$ with components
$
(q_j \cdot p_i \colon Di \to V \to Ej)_{i,j \in \I \times \J}
$.


\begin{Defn}\label{def:2}
  A cocone $p \colon D \to \Delta V$ and a cone $q \colon
  \Delta W \to E$ are said to be \emph{orthogonal}, written $p
  \mathbin \bot
  q$, if for every diagram as in the solid part of
\begin{equation}
\cd[@!@-0.7em]
{
 D \ar[r]^p \ar[d]_h & \Delta V \ar[d]^{k} \ar@{.>}[dl]|{j} \\
 \Delta W \ar[r]_q & E\rlap{ ,}
}\label{eq:2}
\end{equation}
wherein $h$ is a cocone, $k$ is a cone, and  $q \cdot h = k \cdot p
\colon D \rightsquigarrow E$, there exists a unique map $j \colon V
\to W$ as indicated making both triangles commute.
\end{Defn}
Of course, this definition generalises the classical notion of
orthogonality of arrows in a
category~\cite[\S2.1]{Freyd1972Categories}; it also generalises the
notion of orthogonality of \emph{discrete} cones and cocones---ones
indexed by discrete categories---formulated
in~\cite[\S3]{Janelidze2012Weakly}, whose special case dealing with
the orthogonality of an arrow to a discrete cone is already present
in~\cite[\S2.4]{Freyd1972Categories}. 

The orthogonality of arrows underlies the notion of factorisation
system introduced in~\cite[\S2.2]{Freyd1972Categories}; more
generally, the orthogonality of arrows to discrete cones plays a role
in~\cite{Herrlich1974Topological}'s notion of $(E,M)$-category, in
which $E$ is a class of arrows, $M$ an orthogonal class of discrete
cones, and every discrete cone factors as an $E$-map followed by an
$M$-cone. The following definition generalises these notions further
to involve orthogonality of arbitrary small cocones and cones.
\begin{Defn}\label{def:3}
  A \emph{cylinder factorisation system} on a category $\C$ comprises a
  class $\E$ of small cocones---``small'' meaning ``indexed by a small
  category''---and a class $\M$ of small cones, satisfying the
  following properties:
  \begin{enumerate}[(i)]
  \item $\E$ is closed under postcomposition with isomorphisms, and
    $\M$ is closed under precomposition with isomorphisms;
  \item $p \mathbin \bot q$ for all $p \in \E$ and $q \in \M$;
  \item Each small cylinder $r \colon D \rightsquigarrow E$ has a
    factorisation $r = q \cdot p$ with $p \in \E$ and $q \in \M$.
  \end{enumerate}
\end{Defn}
It follows that $\E$ comprises all small cocones $q$ such that $q
\mathbin \bot p$ for all $q \in \M$, and that $\M$ comprises all small
cocones $q$ such that $q \mathbin \bot p$ for all $p \in \E$; and in
fact these two conditions together with (iii) gives an alternate
axiomatisation of cylinder factorisation systems. 
Every cylinder factorisation system $(\E, \M)$ has an underlying
orthogonal factorisation system $(\E_0, \M_0)$---in the sense
of~\cite{Freyd1972Categories}---obtained by restricting to cones and
cocones over diagrams $1 \to \C$. The following result extends one of
the basic facts in that theory to the cylinder setting.

\begin{Lemma}
  \label{lem:3}
  Factorisations in a cylinder factorisation system are essentially
  unique: if the cylinder $r \colon D \rightsquigarrow E$ admits the
  $(\E, \M)$-factorisations $k \cdot p \colon D \to \Delta V \to E$
  and $q \cdot h \colon D \to \Delta W \to E$, then the unique map $j
  \colon V \to W$ as in~\eqref{eq:2} is invertible.
\end{Lemma}
\begin{proof}
  Mirroring~\eqref{eq:2} through the $DE$-axis and applying
  orthogonality again yields a filler $j' \colon W \to V$; now both $j
  \cdot j'$ and $1_W$ fill the square $q \cdot h = q \cdot h$, and so
  must be equal; dually we have $j' \cdot j = 1_V$.
\end{proof}

\begin{Exs}
  \begin{enumerate}[(a)]

  \item If $\C$ is complete, then it admits a cylinder factorisation
    system (small cocones, limit cones). Condition (i) is obvious,
    while (ii) is easy from the universality of a limiting cone. For
    (iii), we may factorise a cylinder $r \colon D \rightsquigarrow E$
    as $p \colon D \to \Delta(\lim E)$ followed by $q \colon
    \Delta(\lim E) \to E$, where $q$ is the limiting cone, and for
    each $i \in \I$, $p_i \colon Di \to \lim E$ is the unique map with
    $q_j \cdot p_i = r_{ij}$ for each $j \in \J$.

  \item Dually, if $\C$ is cocomplete, then it admits a cylinder
    factorisation system (colimit cocones, small cones).

  \item Let $\C$ be complete and cocomplete, and let $(\E_0, \M_0)$ be an
    orthogonal factorisation system on $\C$. We obtain a cylinder
    factorisation system $(\E, \M)$ on $\C$ by
    taking:
    \begin{align*}
       \E &= \{\, p \colon D \to \Delta V \text{ small} :
      \text{the induced } \bar p \colon
      \colim D \to V \text{ is in } \E_0\,\}\\
      \M &= \{\, q \colon \Delta W \to E \text{ small} :
      \text{the induced } \bar q \colon W \to \lim E \text{ is in }
      \M_0\,\}\rlap{ .}
    \end{align*}
    Axiom (i) is clear, while (ii) follows easily on observing that
    diagrams~\eqref{eq:2} correspond bijectively with squares in $\C$
    of the form:
    \begin{equation*}
      \cd[@!@-1.8em] {
        \colim D \ar[r]^-{\bar p} \ar[d]_{\bar h} & V \ar[d]^{\bar k} \ar@{.>}[dl]|{j} \\
        W \ar[r]_-{\bar q} & \lim E\rlap{ .} }
    \end{equation*}
    As for (iii), given $r \colon D \rightsquigarrow E$, we first
    factorise as $q \cdot \ell \colon D \to \Delta(\lim E) \to E$ as
    in (a); then we factorise $\ell$ dually as $f \cdot p \colon D \to
    \Delta(\colim D) \to \Delta(\lim E)$; then we factorise $f = e
    \cdot m \colon \colim D \to V \to \lim E$ with $e \in \E_0$ and $m
    \in \M_0$; and finally take our desired factorisation to be $e
    \cdot p \colon D \to \Delta V$ followed by $q \cdot m \colon
    \Delta V \to E$. It is easy to see that \emph{any} cylinder
    factorisation system on a complete and cocomplete category is
    induced in this way.

  \item Let $\C$ be a complete category which admits (strong epi,
    mono) factorisations and unions of small families of subobjects.
    Call a small cocone $p \colon D \to \Delta V$ \emph{covering} if
    any monomorphism $V' \rightarrowtail V$ through which each $p_i$
    factors is invertible; and call a small cone \emph{monic} if it is
    in $\M$ as defined in (c) for $\M_0$ the class of monomorphisms.
    Now $\C$ admits the cylinder factorisation system (covering
    cocones, mono cones). Axioms (i) and (ii) are straightforward. For
    (iii), given a cylinder $r \colon D \rightsquigarrow E$, we first
    factorise as $q \cdot p \colon D \to \Delta(\lim E) \to E$ as in
    (a); next we (strong epi, mono) factorise each $p_i$ as $m_i \cdot
    e_i \colon Di \twoheadrightarrow Hi \rightarrowtail \lim E$; then
    we form the union $n \colon V \rightarrowtail \lim E$ of the
    subobjects $m_i$ with inclusions $h_i \colon Hi \rightarrowtail
    V$; finally, we obtain our desired factorisation as $h \cdot e
    \colon D \to H \to \Delta V$ followed by $q \cdot n \colon \Delta
    V \to \Delta(\lim E) \to E$. The only non-trivial point is showing
    that $h \cdot e \colon D \to \Delta V$ is covering. So suppose
    that each component $h_i \cdot e_i$ factors through some $g \colon
    V' \rightarrowtail V$. Because each $e_i$ is strongly epic, this
    is equally to say that each $h_i$ factors through $g$; thus each
    $n \cdot h_i = m_i \colon Hi \rightarrowtail \lim E$ factors
    through $n \cdot g \colon V' \rightarrowtail \lim E$; but as $n$
    is the union of the $m_i$'s, $g$ must be invertible as required.

  \item If the \emph{small} category $\C$ bears a cylinder
    factorisation system, then all its $\E$-cocones must be jointly
    epimorphic, and all its $\M$-cones jointly monic, by an adaptation
    of an argument due to Freyd (though see
    also~\cite[Theorem~15.4]{Adamek1990Abstract}). Indeed, suppose
    that $k \colon \Delta V \to E$ is an $\M$-cone, and $f \neq g
    \colon W \to V$ with $k\cdot f = k\cdot g \colon \Delta W \to E$. Let $D$ be
    the discrete diagram comprising $\abs{\mathrm{mor}\ \C}$ copies of
    $W$, let $r \colon D \rightsquigarrow E$ be the cylinder
    comprising $\abs{\mathrm{mor}\ \C}$ copies of the cocone $kf =
    kg$, and let $r = q \cdot p \colon D \rightarrow \Delta U
    \rightarrow E$ be an $(\E, \M)$-factorisation. Then in the diagram
\[
\cd{
  D \ar[d]_\ell \ar[r]^p & 
  \Delta U \ar[d]^q \\ \Delta V \ar[r]_k & E }\] there are at least $2^{\abs{\mathrm{mor}\
   \C}}$ distinct cones $\ell$ yielding commutativity; and so by
orthogonality, at least $2^{\abs{\mathrm{mor}\ \C}}$ distinct maps $U
\to V$ in $\C$, a contradiction.\end{enumerate}\label{ex:1}
\end{Exs}



We now define appropriate notions of morphism between categories
equipped with cylinder factorisation systems. In considering cylinder
factorisation systems on different categories, we will uniformly denote
the classes of cocones and cones by $\E$ and $\M$; normally, context
will make clear which $\E$ and $\M$ are intended, but where confusion
seems possible, we will subscript them with the name of the category
on which they reside.

\begin{Defn}
  \label{def:6}
  We write $\cat{CFS}$ for the $2$-category whose objects are locally
  small categories equipped with a cylinder factorisation system,
  whose $1$-cells are functors $F \colon \C \to \D$ such that $F(\E)
  \subset \E$ and $F(\M) \subset \M$, and whose $2$-cells are
  arbitrary natural transformations. We write $\cat{CFS}_\M$ and
  $\cat{CFS}_\E$ for the corresponding $2$-categories wherein the
  morphisms are required only to preserve $\M$-cones, or only to preserve $\E$-cocones.
\end{Defn}

\begin{Exs}\hfill
  \label{ex:2}
  \begin{enumerate}[(a)]
  \item If the complete $\C$ and $\D$ are equipped with the (all
    cocones, limit cones) cylinder factorisation system, then a
    functor $\C \to \D$ always preserves $\E$-cocones, and preserves
    $\M$-cones precisely when it is continuous. Dually, if the
    cocomplete $\C$ and $\D$ bear the (colimit cocones, all cones)
    cylinder factorisation systems, then a functor between them always
    preserves $\M$-cones and preserves $\E$-cocones just when it is
    cocontinuous. It follows that $\cat{CFS}$ contains as full
    sub-$2$-categories both the $2$-category $\cat{COCTS}$ of
    cocomplete categories and cocontinuous functors, and the
    $2$-category $\cat{CTS}$ of complete categories and continuous functors.
 
  \item If $\C$ and $\D$ are cocomplete, then the condition that a
    morphism $F \colon \C \to \D$ in $\cat{CFS}_\E$ must satisfy can
    be reduced to the requirements that $F(\E_0) \subset \E_0$, and
    that $F$ should preserve colimits ``up to $\E_0$''; meaning that
    each canonical comparison $F \colim D \to \colim FD$ should be in
    $\E_0$. In \cite{Kelly1980A-unified}, Kelly calls this condition
    \emph{preserving the $\E_0$-tightness} of colimit cocones. Of
    course, we have a dual characterisation of morphisms of
    $\cat{CFS}_\M$ between complete categories.

  \item It is easy to see that if $F \dashv G \colon \D \to \C$, and
    $p$ is a cocone in $\C$ and $q$ a cone in $\D$, then $Fp
    \mathbin\bot q$ if and only if $p \mathbin \bot Gq$. It follows
    that, if $\C$ and $\D$ are equipped with cylinder factorisation
    systems, then $F$ preserves $\E$-cocones if and only if $G$
    preserves $\M$-cones. 
  \end{enumerate}
\end{Exs}

We conclude this section with a technical result, necessary in the
sequel, that gives an understanding of the effect of cylinder
factorisation systems on cylinders which, though not small, are
``essentially small'' in a sense now to be described. Recall that a
functor $K \colon \J' \to \J$ is called \emph{initial} if, for each $j
\in \J$, the comma category $K / j$ is connected; which by the
pointwise formula for Kan extensions, is equally to say that the
triangle
\[
\cd{
\J' \ar[rr]^{K} \ar[dr]_{\Delta 1} & \rtwocell{d}{1} & \J
\ar[dl]^{\Delta 1} \\ & \cat{Set}
}
\]
is a left Kan extension. The universal property of Kan extension now
implies that, for each diagram $E \colon \J \to \C$ and $W \in \C$,
precomposition with $K$ induces a bijection
\[
[\J, \cat{Set}](\Delta 1, \C(W, E\thg)) \cong [\J', \cat{Set}](\Delta
1, \C(W, EK\thg))
\]
between cones $q \colon \Delta W \to E$ and cones $qK \colon \Delta W
\to EK$; which in turn implies a bijection between cylinders $r \colon
D \rightsquigarrow E$ and ones $r(1 \times K) \colon D
\rightsquigarrow EK$. Dually, a functor $H \colon \I' \to \I$ is
called \emph{final} if each comma category $i/H$ is connected; which
now implies a bijection between cocones $p \colon D \to \Delta V$ and
ones $pH \colon DH \to \Delta V$, and between cylinders $r \colon D
\rightsquigarrow E$ and ones $r(H \times 1) \colon DH \rightsquigarrow
E$. It immediately follows that:

\begin{Lemma}
  \label{lem:1}
  If $H \colon \I' \to \I$ is final, $K \colon \J' \to \J$ is initial,
  $D \colon \I \to
  \C$ and $E \colon \J \to
  \C$, then for any cocone $p \colon D \to \Delta V$ and any cone $q
  \colon \Delta W \to E$,
  we have $p \mathbin \bot q$ iff $pH \mathbin \bot qK$.
\end{Lemma}

Let us now define a cylinder $r \colon D \rightsquigarrow E$ to be
\emph{essentially small} if the category $\I$ indexing $D$ admits a
final functor from a small category, and the category $\J$ indexing
$E$ admits an initial functor from a small
category. 
In particular, this gives a
notion of essential-smallness for cocones and cones, on identifying
these with degenerate cylinders. 



For the nonce, we will call a structure as in
Definition~\ref{def:3}, but where ``small'' has everywhere been
replaced by ``essentially small'', an \emph{extended cylinder
  factorisation system}. Restricting an extended cylinder
factorisation system to its small cocones and cones yields a cylinder
factorisation system; while in the other direction, we have:

\begin{Prop}
  \label{prop:2}
  Every cylinder factorisation system $(\E, \M)$ on $\C$ is
  the underlying cylinder factorisation system of a unique extended
  cylinder factorisation system $(\overline \E, \overline \M)$;
  moreover, any morphism of cylinder factorisation systems $F \colon
  \C \to \D$ preserves these extended classes, in that $F(\overline \E)
  \subset \overline \E$ and $F(\overline \M) \subset \overline \M$.
\end{Prop}
\begin{proof}
  Given $(\E, \M)$, we define classes of essentially small cocones and
  cones by
  \begin{align*}
    \overline \E &= \{p \colon D \to \Delta V \mid pH 
    \in \E \text{
      for some final } H \colon \I' \to \I \}\\[-2pt]
    \overline \M &= \{q \colon \Delta W \to E \mid qK 
    \in \E \text{ for some initial } K \colon \J' \to \J \}\rlap{ .}
  \end{align*}
  Clearly axiom (i) is satisfied, while (ii) is immediate from
  Lemma~\ref{lem:1}. This same lemma implies that $\overline \E$
  comprises precisely those essentially small cocones orthogonal to
  every cone in $\M$, and vice versa, from which uniqueness of
  $(\overline \E, \overline \M)$ follows easily. The final clause of
  the proposition is immediate from the definitions, and so it remains
  only to show axiom (iii): that each essentially small $r \colon D
  \rightsquigarrow E$ has an $(\overline \E, \overline
  \M)$-factorisation. Given such an $r$, choose a final $H \colon \I'
  \to \I$ and an initial $K \colon \J' \to \J$ with $\I'$ and $\J'$
  small, let $r' = r(H \times K) \colon DH \rightsquigarrow EK$, and
  form $q' \cdot p' \colon DH \to \Delta V \to EK$ an $(\E,
  \M)$-factorisation of the small $r'$. Since $H$ is final and $K$
  initial, there are unique $p \colon D \to \Delta V$ and $q \colon
  \Delta V \to E$ with $pH = p'$ and $qK = q'$, and clearly $p \in
  \overline \E$ and $q \in \overline \M$; finally, since $r(H \times
  K) = r' = q' \cdot p' = qK \cdot pH = (q \cdot p)(H \times K)$, we
  have by finality and initiality of $H$ and $K$ that $r = q \cdot p$, as desired.
\end{proof}
Henceforth, then, there will be no explicit need to speak of extended
cylinder factorisation systems; instead, we modify our notation by
allowing $\E$ and $\M$, which previously denoted the classes of small
cocones and cones of a cylinder factorisation system, to denote
instead the essentially small cocones and cones in the closures
$\overline \E$ and $\overline \M$.


\section{The free cylinder factorisation system}
\label{sec:free-cylind-fact-1}
In this section, we give our first main result, showing that the
Isbell envelope $\Env$ is the free category with a cylinder
factorisation system on $\C$. We begin by constructing the cylinder
factorisation system in question.

\begin{Prop}\label{prop:1}
  For any category $\C$, the Isbell envelope $\Env$ bears a
  cylinder factorisation system whose classes of small cocones and
  cones are given by:
\begin{align*}
\E &= \{\, p \colon D \to \Delta V  \mid \pi_1(p) \text{ is
  colimiting in $\P\C$}\,\}\\
\M &= \{\, q \colon \Delta W \to E  \mid \pi_2(q) \text{ is
  limiting in $\P^\dagger\C$}\,\}\rlap{ ,}
\end{align*}
where $\pi_1 \colon \Env \to \P \C$ and $\pi_2 \colon \Env \to
\P^\dagger \C$ are as in~\eqref{eq:4}.
\end{Prop}

\begin{proof}
  Axiom (i) is clear. For (ii), suppose given a diagram~\eqref{eq:2}
  in $\Env$ with $p \in \E$ and $q \in \M$. Applying $\pi_1$ and
  $\pi_2$ we obtain diagrams
  \begin{equation*}
    \cd[@C-0.8em]{
      D^+ \ar[r]^-{p^+} \ar[d]_{h^+} &
      \Delta(V^+) \ar[d]^{k^+} \ar@{.>}[dl]|{m^+} \\
      \Delta(W^+) \ar[r]_-{q^+} & E^+
    } \qquad \text{and} \qquad
    \cd[@C-0.8em]{
      D^- \ar[r]^-{p^-} \ar[d]_{h^-} &
      \Delta(V^-) \ar[d]^{k^-} \ar@{.>}[dl]|{m^-} \\
      \Delta(W^-) \ar[r]_-{q^-} & E^-
    }
\end{equation*}
in $\P\C$ and in $\P^\dagger\C$ respectively. Now $p^+$ is
colimiting since $p \in \E$; it is thus orthogonal to any small cone,
in particular to $q^+$, and so there is a unique diagonal filler $m^+$
as on the left. Similarly, $q^-$ is limiting since $q \in \M$, whence
there is a unique diagonal filler $m^-$ as on the right. We claim that
$(m^+, m^-) \colon V \to W$ is the required unique diagonal filler in
$\Env$. The only point to check is that each square as on the left
in
\[
\cd[@C+0.5em]{
W^-(b) \times V^+(a) \ar[r]^-{1 \times m^+} \ar[d]_{m^- \times 1} & 
W^-(b) \times W^+(a) \ar[d]^{\xi^W}\\
V^-(b) \times V^+(a) \ar[r]_-{\xi^V} & \C(a,b)
}\qquad
\cd[@C+0.7em]{
W^-(b) \times Di^+(a) \ar[r]^-{1 \times m^+ p^+_i}
\ar[d]_{m^- \times p^+_i} & 
W^-(b) \times W^+(a) \ar[d]^{\xi^W}\\
V^-(b) \times V^+(a) \ar[r]_-{\xi^V} & \C(a,b)
}
\]
commutes. Now, evaluating the colimiting cocone $p^+$ at $a$ yields a
colimiting cocone $(p^+_{i}(a) \colon Di^+(a) \to V^+(a))_{i
  \in \I}$; so by precomposing with these maps, it is enough to show
commutativity of the squares on the right above. But by rewriting the
bottom side using~\eqref{eq:3} for $p^+_i$, this is equally to show
that each square
\[
\cd[@C+1em]{
W^-(b) \times Di^+(a) \ar[r]^-{1 \times m^+ p^+_i}
\ar[d]_{p^-_i m^- \times 1} & 
W^-(b) \times W^+(a) \ar[d]^{\xi^W}\\
Di^-(b) \times Di^+(a) \ar[r]_-{\xi^{Di}} & \C(a,b)
}
\]
commutes, which is so by~\eqref{eq:3} for $h = mp_i$. 

This verifies (ii); and there remains only (iii). Given, then, a
cylinder $r \colon D \rightsquigarrow E$ in $\Env$, we first apply
$\pi_1$ and $\pi_2$ to obtain cylinders $r^+$ and $r^-$ in the
cocomplete $\P\C$ and complete $\P^\dagger\C$, which we then
factor as in (a) and (b) of the preceding section as:
\[
r^+ = 
D^+ \xrightarrow{p^+} \Delta V^+ \xrightarrow{q^+} E^+ \qquad \text{and} \qquad
r^- = 
D^- \xrightarrow{p^-} \Delta V^- \xrightarrow{q^-} E^-
\]
with $p^+$ colimiting and $q^-$ limiting. We next define maps
$\xi^V_{ab} \colon V^-(b) \times V^+(a) \to \C(a,b)$ making $V = (V^+,
V^-, \xi^V)$ into an object of $\Env$. Evaluating the colimiting
$p^+$ and limiting $q^-$ at each object $a$ and $b$ yields colimiting
cocones $(p^+_{i}(a) \colon Di^+(a) \to V^+(a))_{i \in \I}$ and
$(q^-_{j}(b) \colon Ej^-(b) \to V^-(b))_{j \in \J}$ in $\cat{Set}$;
 so to give the $\xi^V_{ab}$'s is equally to give their composites
\[\delta_{abij} \colon Ej^-(b) \times Di^+(a) \to \C(a,b)\]
with the components of these cocones: a family of maps natural in
$a,b,i,j$. To obtain such, consider for each $a,b,i,j$ the
square~\eqref{eq:2} associated to the map $r_{ij} \colon Di \to Ej$ in
$\Env$; the common diagonal of the two sides gives the desired
$\delta_{abij}$'s, whose naturality is  easily checked. The $(\E,
\M)$-factorisation of $r$ in $\Env$ is now given by
\[
D \xrightarrow{(p^+, p^-)} \Delta(V^+, V^-, \xi^V) \xrightarrow{(q^+,
  q^-)} E\rlap{ ;}
\]
the only thing left to check is that the components $(p_i^+, p_i^-)$
and $(q_j^+, q_j^-)$ of the cocone and the cone are in fact maps of
$\Env$. By duality, we need only check the former; thus, that each
square as on the left in
\[
\cd[@C+0.2em]{
V^-(b) \times Di^+(a) \ar[r]^-{1 \times p_i^+} \ar[d]_{p_i^- \times 1} & 
V^-(b) \times V^+(a) \ar[d]^{\xi^V}\\
Di^-(b) \times Di^+(a) \ar[r]_-{\xi^{Di}} & \C(a,b)}
\qquad
\cd[@C+0.2em]{
Ej^-(b) \times Di^+(a) \ar[r]^-{q_j^- \times p_i^+} \ar[d]_{p_i^-q_j^- \times 1} & 
V^-(b) \times V^+(a) \ar[d]^{\xi^{V}}\\
Di^-(b) \times Di^+(a) \ar[r]_-{\xi^{Di}} & \C(a,b)}
\]
commutes. Precomposing with the colimit cocone $(q^-_{j}(b) \colon
Ej^-(b) \to V^-(b))_{j \in \J}$, this is equally to show that each
square as on the right commutes. The upper side is, by definition of
$\xi^V$, the common diagonal of the square~\eqref{eq:2} associated to
$r_{ij}$; but as $p_i^- q_j^- = r_{ij}^-$, the lower side of the above
square is also the lower side of that selfsame~\eqref{eq:2}; whence commutativity.
\end{proof}
We are almost ready to give our first main result. First we need a
preparatory lemma.


\begin{Lemma}
  \label{lem:2}
  For each $X \in \Env$ and $a,b \in \C$, the action of the functors
  $\pi_1$ and $\pi_2$ induce homset isomorphisms
$\pi_1 \colon \Env(Ya, X) \to \P\C(Ya, X^+)$ and $\pi_2 \colon \Env(X,
Yb) \to \P^\dagger\C(X^-, Yb)$.
\end{Lemma}
\begin{proof}
  To give a map $f \colon Ya \to X$ in $\Env$ is to give $f^+ \colon
  \C(\thg, a) \to X^+$ in $\P\C$ together with $f^- \colon \C(a, \thg)
  \to X^-$ in $\P^\dagger\C$ rendering commutative each diagram
  \begin{equation}
\cd[@C+0.2em]{ X^-(b) \times \C(a',a) \ar[r]^-{1 \times f^+}
  \ar[d]_{f^- \times 1} &
  X^-(b) \times X^+(a') \ar[d]^{\xi^X}\\
  \C(a,b) \times \C(a',a) \ar[r]_-{\circ} & \C(a,b)\rlap{ .} }\label{eq:7}
\end{equation}
This forces the components of $f^-$ in $\cat{Set}$ to be given by
$\xi^X_{ab}(\thg, x) \colon X^-(b) \to \C(a,b)$, where $x = f^+(1_a)
\in X^+(a)$. Thus $\pi_1 \colon \Env(Ya, X) \to \P\C(Ya, X^+)$ is
injective; for surjectivity, given any $f^+ \in \P\C(Ya, X^+)$, we may
define $f^-$ in the above manner, and verify naturality and
commutativity in~\eqref{eq:7} using the Yoneda lemma. The case of
$\pi_2$ is dual.
\end{proof}
\begin{Thm}
\label{thm:mainresult}
For any category $\C$, the Yoneda embedding $Y \colon \C \to \Env$
into the Isbell envelope exhibits $\Env$, equipped with the cylinder
factorisation system of Proposition~\ref{prop:1}, as the value at $\C$
of a left biadjoint to the forgetful $2$-functor from $\cat{CFS}$ to
$\cat{CAT}$.
\end{Thm}
\begin{proof}
  We must show that, for any category $\D$ equipped with a cylinder
  factorisation system, the functor
  \begin{equation}
    (\thg) \cdot Y \colon \cat{CFS}(\Env, \D) \to \cat{CAT}(\C, \D)\label{eq:9}
  \end{equation}
  is an equivalence of categories. First we show full fidelity: thus,
  given morphisms of cylinder factorisation systems $F$,~$G \colon
  \Env \to \D$ and a natural transformation $\alpha \colon FY \to GY$,
  we must find a unique $\beta \colon F \to G$ with $\beta Y =
  \alpha$. So given $X \in \Env$, form the category of elements $U
  \colon \el X^+ \to \C$ and dually $V \colon \el X^- \to \C$; by the
  Yoneda lemma, we have a colimit cocone $p^+ \colon YU \to \Delta
  X^+$ in $\P\C$---essentially small as $X^+$ is a small colimit of
  representables---and likewise an essentially small limit cone $q^-
  \colon \Delta X^- \to YV$ in $\P^\dagger \C$. By Lemma~\ref{lem:2},
  these lift to a cocone $p \colon YU \to \Delta X$ and cone $q \colon
  \Delta X \to YV$ in $\Env$, necessarily in $\E$ and $\M$
  respectively.
Now as $F(\E) \subset \E$ and $G(\M) \subset \M$, the diagram
\[
\cd{
FYU \ar[r]^{Fp} \ar[d]_{Gp \cdot \alpha U} & \Delta FX \ar@{.>}[dl]|{\beta_X}
\ar[d]^{\alpha V \cdot
 Fq} \\
\Delta GX \ar[r]_{Gq} & GYV
}
\]
of cocones and cones in $\D$ has top edge in $\E$ and bottom edge in
$\M$. The composites around the two sides agree by naturality of
$\alpha$, and so by orthogonality there is a unique diagonal filler
$\beta_X$ as shown making both triangles commute. If $\beta \colon F
\rightarrow G$ is to extend $\alpha$ and be natural, then it must
render these triangles commutative; so these $\beta_X$'s are the
unique possible choice for an extension, and it remains only to show
their naturality in $X$.

So let $f \colon X \to X'$ in $\Env$; we have the $\E$-cocone $p$ and
$\M$-cone $q$ as before, but now also $p' \colon YU' \to \Delta X'$
and $q' \colon \Delta X' \to YV'$. We also have functors $H = \el f^+
\colon \el X^+ \to \el Y^+$ and $K = \el f^- \colon \el Y^- \to \el
X^-$, satisfying $U'H = U$ and $VK = V'$, and, we claim, rendering
commutative both triangles---and hence the outside---in:
\begin{equation}
\cd{
YU \ar[d]_{p'H} \ar[r]^{p} & \Delta X \ar@{.>}[dl]|f \ar[d]^{qK} \\
\Delta X' \ar[r]_{q'} & YV'\rlap{ .}
}\label{eq:8}
\end{equation}
To see this last claim, note that $\pi_1$ of the top triangle commutes
in $\P\C$ by the Yoneda lemma and definition of $H$, and similarly
$\pi_2$ of the bottom triangle commutes; now apply Lemma~\ref{lem:2}.
Using this, we now show naturality of $\beta$ at $f$; thus that $Gf
\cdot \beta_X = \beta_{X'} \cdot Ff$. By orthogonality it suffices to
show equality after precomposition with the $\E$-cocone $Fp$ and after
postcomposition with the $\M$-cone $Gq'$. For the former, we have
that $Gf \cdot \beta_X \cdot Fp = Gf \cdot Gp \cdot \alpha U = Gp'H
\cdot \alpha U = Gp'H \cdot \alpha U'H = \beta_{X'} \cdot Fp'H =
\beta_{X'} \cdot Ff \cdot Fp$; for the latter, $Gq' \cdot Gf
\cdot \beta_X = GqK \cdot \beta_X = \alpha VK \cdot FqK = \alpha V'
\cdot FqK = \alpha V' \cdot Fq' \cdot Ff = Gq' \cdot \beta_{X'} \cdot
Ff$.

This proves that~\eqref{eq:9} is fully faithful; it remains to show
essential surjectivity. Given $F \colon \C \to \D$, we must exhibit a
map $G \colon \Env \to \D$ of cylinder factorisation systems and a
natural isomorphism $GY \cong F$. For each $X \in \Env$, let $q \cdot
p \colon YU \to \Delta X \to YV$ be its canonical essentially small
cylinder, as above. Since $Y$ is fully faithful, there is a unique
cylinder $r \colon U \rightsquigarrow V$ with $Yr = q \cdot p$; now
let $t \cdot s \colon FU \to \Delta GX \to FV$ be an $(\E,
\M)$-factorisation in $\D$ of the essentially small $Fr \colon FU
\rightsquigarrow FV$. This defines $G$ on objects. On morphisms, let
$f \colon X \to X'$ in $\Env$, and let $p, q, p', q', H$ and $K$ be as
in the preceding paragraph. We have by commutativity in~\eqref{eq:8} and full fidelity of
$Y$ that $r(1 \times K) = r'(H \times 1) \colon U \rightsquigarrow
V'$; whence in the diagram on the left in
\begin{equation*}
\cd{
FU \ar[d]_{s'H} \ar[r]^{s} & \Delta GX \ar@{.>}[dl]|{Gf} \ar[d]^{tK} \\
\Delta GX' \ar[r]_{t'} & FV'
} \qquad \qquad
\cd{
FU \ar[d]_{s''H_gH_f} \ar[r]^{s} & \Delta GX \ar@{.>}[dl]|{G(gf)} \ar[d]^{tK_fK_g} \\
\Delta GX'' \ar[r]_{t''} & FV''
}
\end{equation*}
the composite cylinders $Fr(1 \times K)$ and $Fr'(H \times 1)$ are
equal. Since $s \in \E$ and $t' \in \M$, we induce by orthogonality a
unique filler, as displayed; which gives the action of $G$ on
morphisms. Clearly, when $f = 1_X$, we have $H = K = 1$ and $s = s'$
and $t = t'$ and the unique filler $G1_X$ must be $1_{GX}$. So $G$
preserves identities; as for binary composition, given $f \colon X \to
X'$ and $g \colon X' \to X''$, the map $G(gf)$ is the unique filler
for the square on the right above; but since $Gg \cdot Gf \cdot s = Gg
\cdot s'H_f = s''H_gH_f$ and $t'' \cdot Gg \cdot Gf = t'K_g \cdot Gf =
tK_fK_g$, the map $Gg \cdot Gf$ is also a filler. So $G(gf) = Gg \cdot
Gf$ and $G$ is a functor.

To see that $GY \cong F \colon \C \to \D$, note that the canonical
cylinder $r \colon U \rightsquigarrow V$ in $\C$ associated to $YX \in
\Env$ has $U \colon \C / c \to \C$ and $V \colon c / \C \to \C$ the
forgetful functors from the slice and coslice, and $r_{f:a \to
  c,\, g: c \to b} = gf \colon a \to b$; so in particular,
$r_{1_c,1_c} = 1_c$. Consequently, the chosen factorisation $t \cdot s
\colon FU \to \Delta GYc \to FV$ of $Fr$ in $\D$ involves maps
$s_{1_c} \colon Fc \to GYc$ and $t_{1_c} \colon GYc \to Fc$ with
$t_{1_c} \cdot s_{1_c} = 1_{Fc}$. Now as $1_c$ is terminal in $\C/c$,
the functor $1 \to \C/c$ picking it out is final: whence by
Lemma~\ref{lem:1}, $s_{1_c}$, like $s$, is in $\E$; dually, $t_{1_c}$
is in $\M$. So $t_{1_c} \cdot s_{1_c}$ is an $(\E, \M)$-factorisation
of $1_{Fc}$; but so too is $1_{Fc} \cdot 1_{Fc}$, whence by
Lemma~\ref{lem:3}, $t_{1_c}$ is invertible, and
provides the component
at $c$ of the natural isomorphism $GY \cong F$.

Finally, we must show that $G$ is a map of cylinder factorisation
systems. By duality, we need only show that $G(\E) \subset \E$. So let
$w \colon D \to \Delta X$ be an $\E$-cocone in $\Env$; we must show
that $Gw \colon GD \to \Delta GX$ is an $\E$-cocone in $\D$. Consider
the category $\el D^+$ whose objects are triples $(i \in \I, a \in \C,
d \in Di^+a)$ and whose morphisms $(i,a,d) \to (i',a',d')$ are pairs
of $f \colon i \to i'$ in $\I$ and $k \colon a \to a'$ such that $f
\cdot d = d' \cdot k$. Clearly there is a functor $I \colon \el D^+
\to \I$ sending $(i, a, d)$ to $i$, but there is also a functor $W
\colon \el D^+ \to \el X^+$ sending $(i,a,d)$ to $(a, w_i^+(d))$ and
sending $(f, k)$ to $k$. We claim that \emph{$W$ is
  final}.

Indeed, for any $x \in X^+ a$, the comma category $(a,x) / W$ has
objects being triples of $i \in \I$, $h \colon a \to b$ in $\C$ and $d
\in Di^+ b$ with $x = w_i^+(d) \cdot h$, and morphisms $(i, h, d) \to
(i', h', d')$ being pairs $f \colon i \to i'$ and $k \colon b \to b'$
with $k h = h'$ and $d' \cdot k = f \cdot d$. We must show this
category to be connected. Since any object $(i, h, d)$ admits a map
$(1_i, h)$ from one of the form $(i, 1_a, d')$, it's enough to show
connectedness of the full subcategory on objects of this form. This
subcategory is equally the full subcategory $\A_{a,x} \subset
\el\,(D\thg)^+ a$ on those pairs $(i \in \I, d \in Di^+a)$ with $x =
w_i^+(a)$. Now as $w$ is an $\E$-cocone in $\Env$, its
projection $w^+$ in $\P\C$ is colimiting, which is to say that each
cocone $(w_i^+(a) \colon D_i^+(a) \to X^+a)_{i \in \I}$ is colimiting;
whence $\A_{a,x}$ is connected, $(a,x) / W$ is
connected, and so $W$ is final.

Now, let $\tau \cdot \sigma \colon FU \to \Delta GX \to FV$ be the
factorisation defining $GX$, and for each $i \in \I$, let $t_i \cdot
s_i \colon FU_i \to \Delta GDi \to FV_i$ be the corresponding
factorisation for $GDi$. For each $i \in \I$, let $W_i \colon \el Di^+
\to \el X^+$ be the functor induced by $w_i^+$; note that we have $U_i
= UW_i$ and commuting diagrams of cocones as on the left in
\[
\cd{ FU_i \ar[d]_{s_i} \ar[r]^{\sigma W_i} & \Delta GX \\
\Delta GDi  \ar[ur]_{Gw_i}} \qquad \qquad
\cd{
FUW \ar[r]^{\sigma W} \ar[d]_{s} &  \Delta GX\rlap{ .} \\
GDI \ar[ur]_{GwI}
}
\]
It follows that the natural $s \colon FUW \rightarrow
GDI$ whose component at $(i, a, d) \in \el D^+$ is $(s_i)_{(a,
  d)} \colon Fa \to GDi$ fits into a commuting diagram as on the right
above. We are now ready to prove that $Gw$ is an $\E$-cocone. Suppose given
an $\M$-cone $v$ fitting into a diagram of cocones and cones in $\D$
as on the left in
\[
\cd{
  GD \ar[r]^{Gw} \ar[d]_h & \Delta GX \ar[d]^k  \\
  \Delta W \ar[r]_v & E } \qquad \qquad \cd[@C+0.3em]{ FUW \ar[r]^{GwI \cdot s}
  \ar[d]_{hI \cdot s} & \Delta GX \ar[d]^k
  \ar@{.>}[dl]|m \\
  \Delta W \ar[r]_v & E } \qquad \qquad
\cd[@C+0.3em]{
FU_i \ar[r]^-{Gw_i \cdot s_i} \ar[d]_{h_i \cdot s_i}& \Delta GX \ar@{.>}[dl]|m \ar[d]^k \\
\Delta W \ar[r]_v & E\rlap{ .}
}
\] Whiskering the cocones with $I$ and precomposing with $s$ yields
the commuting diagram in the centre. The top edge therein is $\sigma
W$ which by Lemma~\ref{lem:1} is in $\E$, since $\sigma$ is so and $W$
is final. So by orthogonality there is a unique $m$ as indicated
making both triangles commute. This commutativity is equivalent to
that of the two triangles on the right for every $i \in \I$; wherein
the the condition $m \cdot Gw_i \cdot s_i = h_i \cdot s_i$ for the top
triangle, together with $v \cdot m \cdot Gw_i = k \cdot Gw_i = v \cdot
h_i$, implies that $m \cdot Gw_i = h_i$, since $v \in \M$ and $s_i \in
\E$. So, finally, $m$ is unique such that $v \cdot m = k$ and $m \cdot
Gw = h$, thus a unique filler for the left square, as required.
\end{proof}

\section{Pseudomonadicity}
\label{sec:monadicity}
The preceding result shows that the embedding $Y \colon \C \to \Env$ into the
Isbell envelope is the unit at $\C$ of a biadjunction $\cat{CFS}
\leftrightarrows \cat{CAT}$. This biadjunction induces a pseudomonad
$\En$ on $\cat{CAT}$, and a canonical comparison homomorphism $K
\colon \cat{CFS} \to \En\text-\cat{Alg}$, whose codomain is the
$2$-category of $\En$-pseudoalgebras, algebra pseudomorphisms and
algebra $2$-cells. Recall---for instance, from~\cite[\S
2]{Street1980Fibrations}---that an $\I$-\emph{pseudoalgebra} involves
a morphism $A \colon \Env \to \C$ and invertible $2$-cells $\theta
\colon 1_\C \cong AY$ and $\pi \colon A \cdot \mu_\C \cong A \cdot
\Env[A]$ satisfying two coherence axioms; and that an \emph{algebra
  pseudomorphism} $(\C, A) \to (\D, B)$ involves a morphism $F \colon
\C \to \D$ and an invertible $2$-cell $\varphi \colon B \cdot \Env[F]
\cong F A$, also satisfying two coherence axioms.

Our second main result states that the canonical comparison $K \colon
\cat{CFS} \to \I\text-\cat{Alg}$ is a
biequivalence; in other words, that $\cat{CFS}$ is
\emph{pseudomonadic} over $\cat{CAT}$. We could prove this using the
pseudomonadicity theorem of~\cite{Le-Creurer2002Becks}, but it will be
simpler and more illuminating to construct directly a biequivalence
inverse.

\begin{Thm}
\label{thm:pseudo}
The forgetful $2$-functor $\En\text-\cat{Alg} \to \cat{CAT}$ has a
(strictly commuting) factorisation
\[
\cd{\I\text-\cat{Alg} \ar[rr]^J \ar[dr] & & \cat{CFS} \ar[dl] \\ & \cat{CAT}}
\]
wherein $J$ is a biequivalence $2$-functor satisfying $JK = 1$; it
follows that $K$ is a biequivalence, and so that $\cat{CFS}$ is
pseudomonadic over $\cat{CAT}$.
\end{Thm}
\begin{proof}
  We first introduce some terminology: given a functor $F \colon \C
  \to \D$ and a cylinder factorisation system on $\C$, we say that a
  cocone $p \colon D \to \Delta V$ in $\C$ is \emph{$F$-nearly in}
  $\E$ if, on forming an $(\E, \M)$-factorisation $p = t \cdot s
  \colon D \to \Delta W \to \Delta V$, the map $t$ is inverted by $F$.
  It is easy to see that if $p \in \E$, then $p$ is $F$-nearly in
  $\E$; and that, if $G \colon \B \to \C$ is a map in $\cat{CFS}_\E$,
  then a cocone $p$ in $\B$ is $FG$-nearly in $\E_\B$ iff $Gp$ is
  $F$-nearly in $\E_\C$. Of course, there is the dual notion of a cone
  being $F$-nearly in $\M$, with the corresponding dual results.

  With this in place, we now define $J$ on objects. Let $A \colon \Env
  \to \C$ be an $\En$-pseudoalgebra. We define classes of small cones
  and cocones in $\C$ by:
  \begin{equation}
  \begin{aligned}
    \E &= \{\,p \colon D \to \Delta V \mid Yp \text{ is $A$-nearly in
      $\E_{\Env}$}\,\} \\
    \M &= \{\,q \colon \Delta W \to E \mid Yq \text{ is
      $A$-nearly in $\M_{\Env}$}\,\}\rlap{ ,}
  \end{aligned}\label{eq:6}
\end{equation}
  and claim that this provides the required cylinder factorisation
  system on $\C$. As a first step, we prove that $A \colon \Env \to
  \C$ has $A(\E) \subset \E$ and $A(\M) \subset \M$; by duality we
  need only prove the first. 
  So given $p \in \E_{\I\C}$, we must show that $YAp$ is $A$-nearly in
  $\E_{\I\C}$. By pseudonaturality of the unit of $\En$, we have $YA
  \cong \Env[A] \cdot Y$, so this is equally to show that
  $\Env[A]\cdot Yp$ is $A$-nearly in $\E_{\Env}$. Since $\Env[A]
  \colon \Env[\Env] \to \Env$ is a map of (free) cylinder
  factorisation systems, this is equally to show that $Yp$ is $A \cdot
  \Env[A]$-nearly in $\E_{\Env[\Env]}$; but $A \cdot \Env[A] \cong A
  \cdot \mu_\C$ since $\C$ is a pseudoalgebra, and so this is equally
  to show that $Yp$ is $A \cdot \mu_\C$-nearly in $\E_{\Env[\Env]}$.
  Now as $\mu_\C$ is a map of cylinder factorisation systems, this is
  equally to show that $\mu_\C \cdot Yp$ is $A$-nearly in $\E_{\Env}$;
  finally, since $\mu_\C \cdot Y \cong 1$, this is equally to show
  that $p$ is $A$-nearly in $\E_{\Env}$, which is certainly so if $p
  \in \E_{\I\C}$. 

  We now show that the classes~\eqref{eq:6} verify the axioms
  (i)--(iii) for a cylinder factorisation system on $\C$. (i) is
  trivial; for (iii), given a small cylinder $r \colon D
  \rightsquigarrow E$ in $\C$, we form an $(\E, \M)$-factorisation $Yr
  = q \cdot p$ in $\Env$; by the above, $AYr = Aq \cdot Ap$ is an
  $(\E, \M)$-factorisation in $\C$, and so conjugating by the
  isomorphism $\theta \colon 1_\C \cong AY$ (coming from the
  pseudoalgebra structure of $\C$) we obtain the desired factorisation
  $r = (\theta^{-1} E\cdot Aq) \cdot (Ap \cdot \theta D) \colon D \to
  \Delta V \to E$. It remains to verify (ii). Let $p \in \E_\C$ and $q
  \in \M_\C$ and suppose given a square $q \cdot h = k \cdot p$ as
  in~\eqref{eq:2}. In $\Env$ we may form the diagram on the left
  \[
  \cd{
    YD \ar[r]^{s} \ar[d]_{Yh} & \Delta X \ar[r]^{t} \ar@{.>}[d]^\ell & \Delta YV
    \ar[d]^{Yk} \\
    \Delta YW \ar[r]_{u} & \Delta Y \ar[r]_{v} & YE
  } \qquad \qquad
  \cd{
    AYD \ar[r]^{AYp} \ar[d]_{AYh} & \Delta AYV \ar[d]^{AYk} \ar@{.>}[dl]|{m} \\
    \Delta AYW \ar[r]_{AYq} & AYE\rlap{ ,}
  }
  \]
  wherein both rows are $(\E, \M)$-factorisations and $\ell$ is the
  unique map induced by orthogonality of $s$ and $v$. Since $p \in \E_\C$
  and $q \in \M_\C$, applying $A$
  inverts $u$ and $t$, and so we obtain a diagonal filler for the
  square on the right above by taking $m = (Au)^{-1} \cdot A\ell \cdot
  (At)^{-1}$; conjugating by $\theta \colon 1_\C \cong AY$ now yields
  the required filler $j = \theta^{-1}_W \cdot m \cdot \theta_V \colon
  V \to W$ for the original square~\eqref{eq:2}. To show uniqueness of
  $j$, let $j' \colon V \to W$ be another diagonal filler; then $u
  \cdot Yj' \cdot t \colon \Delta X \to \Delta Y$ fills the rectangle
  on the left above, and so by orthogonality must be $\ell$; whence $A\ell =
  Au \cdot AYj' \cdot At$, so that $m = AYj'$ and so finally $j =
  \theta^{-1}_W \cdot AYj' \cdot \theta_V = j'$.

  This defines $J$ on objects; since $\cat{CFS} \to \cat{CAT}$ is
  faithful on $1$-cells and locally fully faithful, the definition on
  $1$- and $2$-cells is forced, and all that is required is to show
  that any pseudomorphism $F \colon (\C,A) \to (\D,B)$ of
  $\I$-pseudoalgebras preserves the classes of the derived cylinder
  factorisation systems. So let $p$ be a cocone in $\C$ such that $Yp$
  is $A$-nearly in $\E_{\Env}$; we must show that $YFp$ is $B$-nearly
  in $\E_{\Env[\D]}$. By naturality of $Y$, we have $YF \cong \Env[F]
  \cdot Y$, so it's enough to show that $\Env[F] \cdot Yp$ is
  $B$-nearly in $\E_{\Env[\D]}$. Since $\Env[F]$ is a map of cylinder
  factorisation systems, it's enough to show that $Yp$ is $B \cdot
  \Env[F]$-nearly in $\E_{\Env}$; but $B \cdot \I F \cong FA$ as $F$
  is a pseudomorphism, so it's enough to show that $Yp$ is $FA$-nearly
  in $\E_{\Env}$; which is so since $Yp$ is $A$-nearly in $\E_{\Env}$.

  This completes the definition of  $J$; we next show
  that $JK = 1$. This is immediate on $1$- and $2$-cells, since $J$
  and $K$ are both over $\cat{CAT}$ and $\cat{CFS} \to \cat{CAT}$ is
  faithful on $1$- and $2$-cells. To show $JK = 1$ on objects, let
  $\C$ be a category equipped with a cylinder factorisation system;
  then $K\C$ is the pseudoalgebra $A \colon \Env \to \C$ whose
  structure map is obtained by extending the identity $\C \to \C$
  using freeness of $\Env$. Now $JK\C$ is the category $\C$ equipped
  with the cylinder factorisation system $(\E', \M')$ where $\E'$
  comprises those cocones $p$ such that $Yp$ is $A$-nearly in
  $\E_{\Env}$; but as $A \colon \Env \to \C$ is a map of cylinder
  factorisation systems, these are equally the cocones $p$ such that
  $AYp$ is $1_\C$-nearly in $\E$; that is, the $\E$-cocones. Thus $\E'
  = \E$ and similarly $\M' = \M$, so that $JK$ is the identity on
  objects as required.

  Finally, we show that $J$ is a biequivalence. Being a retraction, it
  is clearly surjective on objects; we claim that it also full on
  $1$-cells and locally fully faithful. For the first claim, let $(\C,
  A)$ and $(\D, B)$ be $\En$-pseudoalgebras and $F \colon J(\C, A) \to
  J(\D, B)$ a map of induced cylinder factorisation systems. Then in
  the left square of
  \begin{equation}
\cd{ 
  \Env \ar[r]^-{\Env[F]} \ar[d]_A & \Env[\D] \ar[d]^{B} \\
  \C \ar[r]_F& \D } \qquad \qquad
  \cd[@R+0.4em@C-0.5em]{
    \Env \ar[d]_{A} \ar@/^0.8em/[rr]^-{\Env[F]} 
     & & \Env[\D]
    \ar[d]^{B}\\
    \C \ar@/^0.8em/[rr]^-{F} \dtwocell{rr}{\alpha}
     \ar@/_0.8em/[rr]_-{G}& \dtwocell[0.75]{u}{\varphi_F} &
    \D
  } \ \ = \ \   \cd[@R+0.4em@C-0.5em]{
    \Env \ar[d]_{A} \ar@/^0.8em/[rr]^-{\Env[F]} \dtwocell{rr}{\Env[\alpha]}
    \ar@/_0.8em/[rr]_-{\Env[G]} & & \Env[\D]
    \ar[d]^{B}\\
    \C 
     \ar@/_0.8em/[rr]_-{G}& \dtwocell[0.25]{u}{\varphi_G} &
    \D
  }\label{eq:13}
\end{equation}
all four functors are maps of cylinder factorisation systems.
Moreover, using the unit coherences for $(\C, A)$ and $(\D, B)$ and
pseudonaturality of $Y$, we have an isomorphism $\alpha \colon B\cdot
\Env[F] \cdot Y \cong BYF \cong FAY$, and so, by full fidelity
of~\eqref{eq:9}, a unique invertible $2$-cell $\varphi \colon B \cdot
\Env[F] \cong FA$ with $\varphi Y = \alpha$. This makes $(F, \varphi)
\colon (\C, A) \to (\D, B)$ into an algebra pseudomorphism with $J(F,
\varphi) = F$; the first coherence axiom follows immediately from
$\varphi Y = \alpha$, while the second one, equating two parallel
morphisms in $\cat{CAT}(\Env[\Env], \D)$, follows by fidelity
of~\eqref{eq:9} on observing these morphisms to reside in
$\cat{CFS}(\Env[\Env], \D)$ and to have the same precomposite with $Y
\colon \Env \to \Env[\Env]$. It remains to show local full fidelity of
$J$; thus, that for any pair of algebra pseudomorphisms $(F,
\varphi_F), (G, \varphi_G) \colon (\C, A) \to (\D, B)$ and any
$2$-cell $\alpha \colon F \Rightarrow G$, the pasting equality above
right holds. This follows, again, by observing these pastings to
describe parallel morphisms in $\cat{CFS}(\Env, \D)$ which coincide on
precomposition with $Y \colon \C \to \Env$.
\end{proof}

\newcommand{\Lax}{\En\text-\cat{Alg}_{\ell}}
\newcommand{\Colax}{\En\text-\cat{Alg}_{c}}
\section{Lax and colax morphisms}
\label{sec:2-categ-aspects}
As well as the $2$-category $\En\text-\cat{Alg}$, we also have the
larger $2$-categories $\Lax$ and $\Colax$ whose objects are again
pseudoalgebras, but whose $1$- and $2$-cells are now the \emph{lax} or
\emph{colax} algebra morphisms and the algebra $2$-cells between them.
A lax algebra morphism $(\C,A) \rightsquigarrow (\D, B)$ comprises a
functor $F \colon \C \to \D$ and a potentially non-invertible $2$-cell
$\varphi \colon B \cdot \Env[F] \Rightarrow FA$ satisfying two
coherence axioms; a \emph{colax} morphism is similar, but with the
orientation of the non-invertible $\varphi$ now reversed. Our final
result identifies the lax and colax $\En$-algebra morphisms as the
functors preserving only $\M$-cones and only $\E$-cocones
respectively. As in the preceding section, we could proceed by
applying a general theorem, in this case the two-dimensional
monadicity theorem of~\cite{Bourke2014Two-dimensional}; but as there,
it will be simpler and more illuminating to give the constructions
directly.

\begin{Thm}\label{thm:3}
  The factorisation of $\En\text-\cat{Alg} \to \cat{CAT}$ through
  $\cat{CFS}$ extends to a factorisation of $\Lax \to \cat{CAT}$
  through $\cat{CFS}_\M$ and to one of $\Colax \to \cat{CAT}$
  through $\cat{CFS}_\E$:
\[
\cd{
  \I\text-\cat{Alg} \ar[rr]^J \ar[d] \ar@{-->}[ddr]  & & \cat{CFS} \ar[d] \ar@{-->}[ddl]\\ 
  \Lax \ar[rr]^{J_\ell} \ar[dr] & & \cat{CFS}_\M\ar[dl] \\ 
& \cat{CAT}}\qquad \qquad \qquad
\cd{
  \I\text-\cat{Alg} \ar[rr]^J \ar[d] \ar@{-->}[ddr]  & & \cat{CFS} \ar[d] \ar@{-->}[ddl]\\ 
  \Colax \ar[rr]^{J_c} \ar[dr] & & \cat{CFS}_\E\ar[dl] \\ 
& \cat{CAT}}
\]
wherein $J_\ell$ and $J_c$ are biequivalences.
\end{Thm}
\begin{proof}
  By duality, we consider only the lax case. First we extend $J$ to
  $J_\ell$; of course, $J_\ell$ must agree with $J$ on objects, and as
  before the definition is forced on $1$- and $2$-cells; so the only
  work is showing that, if $(F, \varphi) \colon (\C, A)
  \rightsquigarrow (\D, B)$ is a lax algebra map, then $F$ sends
  $\M_\C$-cones to $\M_\D$-cones. 
  So let $q \colon \Delta V \to E$ be an $\M_\C$-cone; we must show
  $Fq \in \M_\D$. Let $Yq = t \cdot s \colon \Delta YV \to \Delta W
  \to YE$ be an $(\E, \M)$-factorisation in $\Env$, and consider the
  commuting diagram on the left in
\[
\cd[@C+0.2em]{
\Delta B(\Env[F])YV \ar[d]_{B(\Env[F])s} \ar[r]^-{\varphi_{YV}} &
\Delta FAYV \ar[d]^{FAs}  \\
\Delta B(\Env[F])W \ar[d]_{B(\Env[F])t} \ar[r]^-{\varphi_{W}} &
\Delta FAW \ar[d]^{FAt} \\
B(\Env[F])YE \ar[r]^-{\varphi YE} &
FAYE
} \quad
\cd[@C+2em]{
\Delta FAYV \ar[d]|{FAYq} \ar@{ >->}[r]^-{B(\Env[F])s \cdot \varphi_{YV}^{-1}} &
\Delta B(\Env[F])W \ar[d]^{B(\Env[F])t} \ar@{->>}[r]^-{FAs^{-1} \cdot \varphi_{W}} & 
\Delta FAYV \ar[d]|{FAYq}\\
FAYE \ar@{ >->}[r]^-{\varphi YE^{-1}} &
B(\Env[F])YE \ar@{->>}[r]^-{\varphi YE} &
FAYE\rlap{ .}
}
\]
To say $q \in \M_\C$ is to say that $Yq$ is $A$-nearly in $\M_{\Env}$:
so $FAs$ is invertible, and by the unit coherence axiom for a lax
morphism so too are $\varphi_{YV}$ and $\varphi YE$. Moreover
$B(\Env[F])t \in \M_{\D}$ since $t \in \M_{\Env}$ and $B \cdot
\Env[F]$ is a map of cylinder factorisation systems. Thus the diagram
on the right 
exhibits $FAYq$ as being a retract of an $\M_\D$-cone and
so, by an easy argument, itself an $\M_\D$-cone; 
finally, since $AY \cong 1$, we have $Fq \cong FAYq$ an $\M_\D$-cone as
required.

This completes the definition of $\J_\ell$, and it remains to show
that it is a biequivalence. Of course, it is surjective on objects,
since $J$ is; we claim it is also full on $1$-cells and locally fully
faithful. We use the fact---generalising full fidelity
of~\eqref{eq:9}---that for any $F \in \cat{CFS}_\E(\Env,\D)$ and $G
\in \cat{CFS}_\M(\Env, \D)$, the function
\begin{equation}
Y \cdot (\thg) \colon \cat{CFS}_\M(\I\C, \D)(F, G) \to \cat{CAT}(\C, \D)(F, G)\label{eq:14}
\end{equation}
is invertible; the proof is precisely the first two paragraphs of the
proof of Theorem~\ref{thm:mainresult}, noting that there we only
needed that $F(\E) \subset \E$ and that $G(\M) \subset \M$. To show
$J_\ell$ is full on $1$-cells, let $(\C, A)$ and $(\D, B)$ be
$\En$-pseudoalgebras and let $F \colon J_\ell(\C,A) \to J_\ell(\D, B)$
in $\cat{CFS}_\M$; then in the left square of~\eqref{eq:13}, the maps
along the upper side are in $\cat{CFS}$, and those along the lower
side in $\cat{CFS}_\M$; so by invertibility of~\eqref{eq:14}, the
isomorphism $\alpha \colon B\cdot \Env[F] \cdot Y \cong BYF \cong FAY$
induces a unique $2$-cell $\varphi \colon B \cdot \Env[F] \Rightarrow
FA$ with $\varphi Y = \alpha$. Using injectivity of~\eqref{eq:14} and
arguing as in the final paragraph of Theorem~\ref{thm:pseudo}, we may
show that this makes $(F, \varphi) \colon (\C, A) \rightsquigarrow
(\D, B)$ into a lax algebra morphism with $J_\ell(F, \varphi) = F$; so
$J_\ell$ is full on $1$-cells. In a similar manner, the argument
showing local full fidelity of $J$ generalises using~\eqref{eq:14} to
one showing local full fidelity of $J_\ell$.
\end{proof}

\section{$(\Phi, \Psi)$-cylinder factorisation systems}
\label{sec:generalisations}
The definition of cylinder factorisation system involves
factorisations for all \emph{small} cylinders---ones indexed by small
categories. However, we could equally well have required
factorisations only for finite cylinders, say, or only for discrete
ones. In this final section, we exhibit such variant notions as the
pseudoalgebras for corresponding variants of the Isbell monad,
obtained by replacing the pseudomonads $\P$ and $\P^\dagger$ used in
its construction by suitable \emph{full submonads} thereof.

By a full submonad $\Ss$ of a pseudomonad $\T$ on $\cat{CAT}$, we mean
the choice, for each category $\C$, of a full subcategory $\Ss\C
\subset \T\C$, with these choices being closed under the pseudomonad
structure of $\T$ in an obvious sense. In the case of $\P$ and
$\P^\dagger$, full submonads $\Phi \subset \P$ and $\Psi \subset
\P^\dagger$ correspond to \emph{saturated classes of weights} for
colimits or limits in the sense of~\cite{Albert1988The-closure} (there
called \emph{closed} classes); the corresponding $\Phi$- or
$\Psi$-pseudoalgebras are categories admitting all $\Phi$-weighted
colimits or all $\Psi$-weighted limits, respectively. Relative to a
choice of full submonads $\Phi \subset \P$ and $\Psi \subset
\P^\dagger$, we may construct a modified Isbell envelope whose value
at a category $\C$ is obtained as a pullback
\begin{equation}
\cd{
\I_{\Phi, \Psi}(\C) \ar[r] \ar[d] & \Env \ar[d]^{(\pi_1,\pi_2)} \\
\Phi\C \times \Psi\C \ar[r] & \P\C \times \P^\dagger \C\rlap{ .}
}\label{eq:16}
\end{equation}
Note that each $\I_{\Phi, \Psi}(\C) \to \I\C$ may be taken to be the
inclusion of a full subcategory; if we do so, then it is easy to see
that these full inclusions assemble together to yield a full submonad
$\I_{\Phi, \Psi} \subset \I$---whose pseudoalgebras we now
characterise.

A diagram $D \colon \I \to \C$ will be called a \emph{$\Phi$-diagram}
if it admits a factorisation as on the left below for some $\varphi
\in \Phi\C$. Dually, $E \colon \J \to \C$ is a \emph{$\Psi$-diagram}
if for some $\psi \in \Psi\C$ it admits a factorisation as on the
right:
\begin{equation}
  \cd[@!C]{ \I \ar[rr]^{H\ \text{final}} \ar[dr]_{D} & & \el \varphi \ar[dl]^{\pi} \\ &
  \C } \qquad \qquad \cd[@!C]{ \J \ar[rr]^{K\ \text{initial}}
  \ar[dr]_{E} & & \el \psi\rlap{ .}
  \ar[dl]^{\pi} \\ & \C }\label{eq:17}
\end{equation}
A \emph{$(\Phi, \Psi)$-cylinder factorisation system} is now defined
identically to a cylinder factorisation system, except that the cones,
cocones and cylinders appearing in the definition are restricted to
those whose domains and codomains are $\Phi$- and $\Psi$-diagrams
respectively. Categories equipped with $(\Phi, \Psi)$-cylinder
factorisation systems are the objects of a $2$-category
$\cat{CFS}_{\Phi, \Psi}$, whose maps are, as before, functors
preserving the cocones and cones.

The proof of the following result follows precisely the arguments of
the preceding sections, but with $\Phi$ and $\Psi$ everywhere
replacing $\P$ and $\P^\dagger$, and with $\Phi$-weighted colimits and
$\Psi$-weighted limits replacing arbitrary colimits and limits. There
is also an analogue of Theorem~\ref{thm:3}, which we do not trouble to
state, characterising the lax and oplax algebra morphisms in terms of
maps preserving only cones or only cocones.

\begin{Thm}
  \label{thm:4}
  Given full submonads $\Phi \subset \P$ and $\Psi \subset
  \P^\dagger$, we have a pseudomonadic adjunction
\begin{equation*}
\cd[@C+2.5em]{
\cat{CFS}_{\Phi, \Psi} \ar@<-1pt>[r]_{} \ar@<4.5pt>@{}[r]|{\bot} \ar@<0pt>@/^0.8em/@{<-}[r]^{} &
\cat{CAT}
}\
\end{equation*}
whose unit at $\C$ may be taken to be the restricted Yoneda embedding
$Y \colon \C \to \I_{\Phi, \Psi}(\C)$.  
\end{Thm}

In practice, the notions of $\Phi$-diagram and $\Psi$-diagram tend to
encompass slightly more than we would intuitively expect. For example,
when $\Phi = 1_\cat{CAT}$, the $\Phi$-diagrams are those $D \colon \I
\to \C$ which admit an absolute colimit in $\C$, rather than simply
those $D \colon 1 \to \C$ indexed by the terminal category. Towards
rectifying this, we define a class $\A$ of $\Phi$-diagrams to be
\emph{generating} if, for every $\varphi \in \Phi\C$, there is some $D
\in \A$ fitting into a diagram as to the left of~\eqref{eq:17}; we
define a generating class $\B$ of $\Psi$-diagrams dually. If $\A$ and
$\B$ are generating classes, then by using Lemma~\ref{lem:1} and
arguing as in Proposition~\ref{prop:2}, we may show that a $(\Phi,
\Psi)$-cylinder factorisation system is completely and uniquely
determined by its cocones, cones, and cylinder factorisations with
respect to diagrams in $\A$ and $\B$.

\begin{Exs}
  \label{ex:3}
  \begin{enumerate}[(a)]
  \item Let $\Phi = \P$ and let $\Psi = \F$ be the pseudomonad for
    \emph{finite} limits---for which $\F\C$ is given by the closure of
    the representables under finite limits in $\P^\dagger \C$---with
    as generating class of $\Psi$-diagrams all diagrams indexed by a
    finite category. In this case, a $(\Phi, \Psi)$-cylinder
    factorisation system involves factorisations for all cylinders
    with \emph{finite} codomain. For example, any regular category
    with pullback-stable unions of subobjects admits a $(\Phi,
    \Psi)$-cylinder factorisation system given by (covering cocones,
    jointly monic cones).

  \item Let $\Phi = \Psi = 1_{\cat{CAT}}$, and take as generating
    classes of $\Phi$- and $\Psi$-diagrams just those indexed by the
    terminal category $1$. Then a $(\Phi, \Psi)$-cylinder
    factorisation system is precisely an orthogonal factorisation
    system; moreover, $\I_{\Phi, \Psi}(\C)$ is the arrow category
    $\C^\cat 2$, and a short calculation shows the pseudomonad
    structure of $\I_{\Phi, \Psi}$ to be that of the ``squaring''
    monad $(\thg)^\cat 2$ of~\cite{Korostenski1993Factorization}.
    Thus we reconstruct the main result of \emph{ibid}., identifying
    orthogonal factorisation systems with $(\thg)^\cat 2$-pseudoalgebras.
  \item Let $\Phi = \cat{Fam}_\Sigma$ and $\Psi = \cat{Fam}_\Pi$ be
    the pseudomonads whose components at $\C$ comprise the coproducts,
    respectively products, of representables in $\P\C$ and $\P^\dagger
    \C$, and take as generating classes of $\Phi$- and $\Psi$-diagrams
    just those indexed by discrete categories. A $(\Phi,
    \Psi)$-cylinder factorisation system now involves factorisations
    of small discrete cylinders---\emph{arrays} in the terminology
    of~\cite{Shulman2012Exact}---into discrete cones and discrete
    cocones, and the notion of orthogonality involved is precisely
    that of~\cite[\S 3]{Janelidze2012Weakly}. In this case, the fact
    that $\I_{\Phi, \Psi}(\C)$ is the free $(\Phi, \Psi)$-cylinder
    factorisation system is quite palpable, since its objects
    are \emph{precisely} the small discrete cylinders in $\C$.
    
    
  \item Let $\Phi = 1_\cat{CAT}$ and $\Psi = \cat{Fam}_\Pi$, with
    generating classes of $\Phi$- and $\Psi$-diagrams as before. In
    this case, a $(\Phi, \Psi)$-cylinder factorisation system involves
    factorisations of discrete cones into $\E$-maps followed by
    $\M$-cones; it is thus a \emph{factorisation structure for small
      sources} in the sense
    of~\cite[Exercise~15J]{Adamek1990Abstract}. As in the preceding
    example, $\I_{\Phi, \Psi}(\C)$ has a simple description as the
    category of all small discrete cones in $\C$.

  \item Let $\Phi = (\thg)_\bot$ be the pseudomonad which freely
    adjoins an initial object, with as generating class of
    $\Phi$-diagrams precisely those indexed by $0$ or $1$; and let
    $\Psi = 1_\cat{CAT}$, with generating class as before. In this
    case, a $(\Phi, \Psi)$-cylinder factorisation system is an
    orthogonal factorisation system in which, additionally, every
    object admits an $\M$-map from an object orthogonal to every
    $\M$-map. As in Examples~\ref{ex:1}(c), this second condition
    follows automatically from the first in the presence of an initial
    object; but there are important cases where initial objects do not
    exist. For example, a category $\C$ admits a $(\Phi,
    \Psi)$-cylinder factorisation system with $\M$ the class of
    \emph{all} maps just when every $A \in \C$ admits a map from a
    \emph{strict generic}~\cite{Weber2004Generic}---an object $G$ such that, for every $X \in
    \C$, the action of $\mathrm{Aut}(G)$ on $\C(G,X)$ is free and
    transitive. 
  \end{enumerate}
\end{Exs}

\bibliography{bibdata}

\end{document}